\documentclass{article}
 \usepackage{amsmath}

\newtheorem{thm}{Theorem}[section]
\newtheorem{cor}[thm]{Corollary}

\newtheorem{lem}[thm]{Lemma}

\newtheorem{rem}[thm]{Remark}

\newcommand{\be}{\begin{equation}}
\newcommand{\ee}{\end{equation}}
\newcommand{\ben}{\begin{enumerate}}
\newcommand{\een}{\end{enumerate}}
\newcommand{\beq}{\begin{eqnarray}}
\newcommand{\eeq}{\end{eqnarray}}
\newcommand{\beqn}{\begin{eqnarray*}}
\newcommand{\eeqn}{\end{eqnarray*}}

\newcommand{\pa}{\partial}

\newcommand{\qed}{\hspace*{\fill}Q.E.D.}  

\begin{document}
\title{On a Class of Finsler Metrics of Scalar Flag Curvature }
\author{Guojun Yang}
\date{}
\maketitle
\begin{abstract}

We have shown that the Beltrami Theorem in Riemannian geometry is
still true for square metrics if the dimension $n\ge 3$, namely,
an $n(\ge 3)$-dimensional square metric is locally projectively
flat if and only if it is of scalar flag curvature. In this paper,
we go on with the study of  the Beltrami Theorem for a larger
class of $(\alpha,\beta)$-metrics $F=\alpha\phi(\beta/\alpha)$
including square metrics, where $\phi(s)$ is determined by a
family of known ODEs satisfied by projectively flat
$(\alpha,\beta)$-metrics. For this class, we prove that the
Beltrami Theorem  holds if $\beta$ is closed, and in particular,
we prove that $\beta$ must  be closed for a subclass with
$\phi(s)$ being a polynomial of degree two. Further,  we obtain
the local and in part the global classifications to those metrics
of scalar flag curvature.


\end{abstract}

\section{Introduction}

It is the Hilbert's Fourth Problem to study and classify locally
projectively flat metrics. The Beltrami Theorem in Riemannian
geometry states that a Riemannian metric is locally projectively
flat if and only if it is of constant sectional curvature. For
Finsler metrics,  the flag curvature is a natural extension of the
sectional curvature in Riemannian geometry, and  every
two-dimensional Finsler metric must be of scalar flag curvature.
It is known that every locally projectively flat Finsler metric is
of scalar flag curvature. However, the converse is not true. For
example, Randers metrics of constant flag curvature  are not
necessarily locally projectively flat (\cite{BRS}).  Therefore, it
is a natural problem to study and classify Finsler metrics of
scalar flag curvature. This problem is far from being solved for
general Finsler metrics. Thus we shall investigate some special
classes of Finsler metrics. Recent studies on this problem are
concentrated on Randers metrics and square metrics.

Randers metrics are among the simplest Finsler metrics in the
following form
$$ F= \alpha+\beta, $$
 where $\alpha$ is a
Riemannian metric  and
 $\beta$ is a 1-form satisfying $\|\beta\|_{\alpha}<1$.
After many mathematician's efforts (\cite{BR1}
\cite{Ma2}--\cite{Shen2} \cite{YS}), Bao-Robles-Shen finally
classify Randers metrics of constant flag curvature  by using the
navigation method (\cite{BRS}). Further, Shen-Yildirim classify
Randers metrics of weakly isotropic flag curvature (\cite{SYi}).
There  are Randers metrics
 of scalar flag curvature which are neither of weakly isotropic flag curvature nor  locally
 projectively flat (\cite{CZ} \cite{SX}).   So far,  the problem of classifying Randers metrics of scalar flag
curvature still remains open.

Recently, square metrics have been shown to have many special
geometric properties. A square metric is defined in the following
form
\[  F =\frac{ (\alpha+\beta)^2}{\alpha},\]
where $\alpha$ is a Riemannian metric and $\beta$ is a $1$-form
with $\|\beta\|_{\alpha}<1$. In \cite{SY}, Shen-Yildirim determine
the local structure of all locally projectively flat square
metrics of constant flag curvature.  L. Zhou shows  that a square
metric of constant flag curvature must be locally projectively
flat (\cite{Z}). Later on, the present author and Z. Shen further
prove that a square metric in dimension $n\ge 3$ is of scalar flag
curvature if and only if it is locally projectively flat, and they
also classify  closed manifolds with a square metric of scalar
flag curvature in dimension $n\ge 3$ (\cite{SYa}).

In this paper, we extend the class of square metrics to a larger
class which is determined by the ODE (\ref{y1}). Then we show the
local and global structures of those metrics which are of scalar
flag curvature in dimension $n\ge 3$.

Now we  state our main results as follows.

\begin{thm}\label{th01}
  Let $F=\alpha\phi(\beta/\alpha)$ be an
  $(\alpha,\beta)$-metric on an  $n(\ge 3)$-dimensional
 manifold $M$, where $\alpha =\sqrt{a_{ij}(x)y^iy^j}$ is Riemannian and $\beta = b_i(x) y^i$ is a
 1-form, and $\phi(s)$ satisfies $\phi(0)=1$ and the following ODE
 \be\label{y1}
    \big\{1+(k_1+k_3)s^2+k_2s^4\big\}\phi''(s)=(k_1+k_2s^2)\big\{\phi(s)-s\phi'(s)\big\},
   \ee
where $k_1,k_2,k_3$ are constant with $k_2\ne k_1k_3$. Assume
$\beta$ is closed. Then $F$
    is of scalar flag curvature if and only if $F$ is locally
    projectively flat.
\end{thm}

 The ODE (\ref{y1}) appears in characterizing an
 $(\alpha,\beta)$-metric which is Douglasian or locally projectively
 flat (\cite{LSS} \cite{Shen1} \cite{Y1}).   The scalar flag curvature {\bf K} in Theorem \ref{th01} can
 be determined (see (\ref{y6}) below).

 If $k_2=k_1k_3$, then $F$ is of Randers
 type. We call an $(\alpha,\beta)$-metric $F=\alpha\phi(s)$ of Randers type,
  if  $\phi(s)= a_1s+\sqrt{1+ks^2}$ for some constants $a_1,k$.
 A metric  of Randers
 type is
 essentially a Randers metric if $k> -1/b^2$. As shown above, it is still an open
 problem to classify Randers metrics of scalar flag curvature.

When $k_1=2,k_2=0,k_3=-3$ and $\phi'(0)=2$, we get
$F=(\alpha+\beta)^2/\alpha$, which is called a square metric. In
this case, the conclusion in Theorem \ref{th01} is still true
without  the condition that $\beta$ is closed (see \cite{SYa},
also see Theorem \ref{th3} below), and by (\ref{y6}) below we have
 \be\label{y008}
{\bf K}=\frac{\alpha}{F^2}\Big\{ [ \lambda + \tau^2 (5+4 b^2)]
\alpha +(\frac{\eta}{2}-3\tau^2)\beta\Big\},
 \ee
which is also obtained in \cite{SYa}. Theorem \ref{th01} might be
also true without the condition that $\beta$ is closed, but it
seems hard to be proved.

 To prove Theorem
\ref{th01}, we first characterize  the metric $F$ in Theorem
\ref{th01} which are of scalar flag curvature in terms of the
covariant derivatives $b_{i|j}$ and the Riemann curvature
$\bar{R}^i_{\ k}$ of $\alpha$ (see Theorem \ref{th1} below). Then
Theorem \ref{th01} follows directly from Theorem \ref{th1}. We
also show that the inverse of Theorem \ref{th1} is true (see
Theorem \ref{th31} below). Besides, in Section \ref{sec5} below,
we use Theorem \ref{th1} to give the classification for an
$(\alpha,\beta)$-metric $F=\alpha\phi(\beta/\alpha)$ with
$\phi(s)$ satisfying (\ref{y1}) which is locally projectively flat
with constant flag curvature (see Corollary \ref{cor51} below)
(cf. \cite{LS} \cite{Y2}).

More important is that, based on Theorem \ref{th1}, we can use the
deformation determined by (\ref{gjcw58})--(\ref{gjcw61}) to obtain
some rigidity results (see Theorem \ref{th2} and Theorem \ref{th3}
below). We will prove in Section \ref{sec4} that, under the
deformation determined by (\ref{gjcw58})--(\ref{gjcw61}) (cf.
\cite{Yu} \cite{Yu1}), if $\alpha$ and $\beta$ satisfy
(\ref{y2})--(\ref{y5}) below, then $h=h_{\mu}$ is a Riemann metric
of constant sectional curvature (put as $\mu$) and $\rho$ is a
closed 1-form which is conformal with respect to $h_{\mu}$.
 Put $h_{\mu}=\sqrt{h_{ij}y^iy^j}$, and
the covariant
    derivatives $p_{i|j}$ of $\rho=p_iy^i$ with
    respect to $h_{\mu}$ satisfy
 \be\label{y9}
 p_{i|j}=-2ch_{ij}
 \ee
 for some scalar function $c=c(x)$ (see (\ref{gjcw63}) below). Let $\nabla c$ be
 the gradient of $c$ with respect to $h_{\mu}$, and then
  \be\label{y010}
 \delta:=\sqrt{\|\nabla c\|^2_{h_{\mu}}+\mu c^2}, \ \ (\mu>0),
  \ee
  is a constant (Lemma \ref{lem61} below). We need $c$ and
  $\delta$ in the following two theorems.

Now if we assume the manifold $M$ is compact without boundary,
then we have the following rigidity theorem based on Theorem
\ref{th1} and the deformation (\ref{gjcw61}).

\begin{thm}\label{th2}
Let $F=\alpha\phi(\beta/\alpha)$ be an
  $(\alpha,\beta)$-metric on an  $n(\ge 3)$-dimensional
 compact  manifold $M$ without boundary, where $\phi(s)$ satisfies
 (\ref{y1}) with $\phi(0)=1$ and $k_2\ne k_1k_3$.

\ben
    \item[{\rm (i)}]  Suppose $F$ is of constant flag
    curvature and $\beta$ is closed. Assume
    \be\label{y0010}
    2+(2k_1-a_1^2)b^2>0, \ \ (a_1:=\phi'(0)),
    \ee
    on the whole $M$ if $F$ is essentially a square metric. Then $F=\alpha$ is Riemannian, or $F$ is
       locally Minkowskian. In the latter case, $F$ is
       flat-parallel ($\alpha$ is flat and $\beta$ is parallel
       with respect to $\alpha$).

    \item[{\rm (ii)}]  Suppose $F$ is of scalar flag
    curvature and $\beta$ is closed (or equivalently, $F$ is locally projectively flat). Then one of the following cases holds:
     \ben
       \item[{\rm (iia)}] If $\mu<0$ for $h$ in (\ref{y9}), then $F=\alpha \ (=h_{\mu})$ is Riemannian.

       \item[{\rm (iib)}] If $\mu=0$ for $h$ in (\ref{y9}), then $F$ is  flat-parallel.

      \item[{\rm (iic)}] If $\mu>0$ for $h$ in (\ref{y9}), then $\alpha$
      and $\beta$ can be written as
       \be\label{ycw11}
     \alpha^2=u^{-1}\big[h_{\mu}^2-4(\mu w)^{-2}vc_0^2\big],\ \
     \beta=2(\mu w)^{-1} c_0,
       \ee
       where $c$ is given by  (\ref{y9}),
       and $c_i:=c_{x^i},c_0:=c_iy^i$, and $u=u(b^2),v=v(b^2),w=w(b^2)$ are determined
      properly by (\ref{gjcw58})--(\ref{gjcw60}) below (for example, by (\ref{gjcw071}) below).
        Further, the scalar flag curvature {\bf K} satisfies
        (\ref{y6}) below with
      \be\label{ycw12}
     \tau=-2uw^{-1}c, \ \
     \lambda=4w^{-2}u\big\{v-k_1(2+k_2b^2)u\big\}c^2+\mu u.
      \ee
     \een
\een
\end{thm}

 It seems difficult to evaluate the lower
and upper bounds of the {\bf K} in Theorem \ref{th2}(iic). But in
a special case for a square metric $F=(\alpha+\beta)^2/\alpha$, it
is shown in \cite{SYa} that {\bf K} satisfies

\be\label{ycw13}
      16\ {\bf
       K}=(\delta^2\mu^{-1}+\mu/4)\mu^3\Big[(1+s)(\delta^2\mu^{-1}+\mu/4-c^2)\Big]^{-3},
         \ee
         \be\label{ycw14}
     \frac{\big(\sqrt{4\delta^2+\mu^2}-2\delta\big)^3}{\mu\sqrt{4\delta^2+\mu^2}}
         \le{\bf  K}\le
         \frac{\big(\sqrt{4\delta^2+\mu^2}+2\delta\big)^3}{\mu\sqrt{4\delta^2+\mu^2}}.
         \ee

\bigskip

Similarly, Theorem \ref{th2} might be also true without the
condition that $\beta$ is closed. For a special case that
$\phi(s)=1+a_1s+\epsilon s^2$ is a polynomial of degree two, we
can prove that $\beta$ is closed, which is shown in the following
theorem (we may put $\epsilon=\pm 1$).

\begin{thm}\label{th3}
 Let $F=\alpha+a_1\beta/\alpha+\epsilon \beta^2/\alpha$ be an
 $(\alpha,\beta)$-metric on an  $n$-dimensional
 manifold $M$, where $a_1$ is a constant and $\epsilon=\pm
 1$.
 \ben
    \item[{\rm (i)}] For $n\ge 3$, $F$ is of scalar flag curvature if and only if $F$ is
 locally projectively flat. In this case,
 (\ref{y2})--(\ref{y007}) below
 hold with $k_1=2\epsilon,k_2=0 ,k_3=-3\epsilon$.

\item[{\rm (ii)}] Let $M$ be
 compact  without boundary. If $F$ is of constant flag curvature, then
 the conclusions in Theorem \ref{th2}(i) hold for $n\ge 2$; and if $F$ is of scalar flag curvature, then
  the conclusions in Theorem \ref{th2}(ii)
 hold for $n\ge 3$. Further for a suitable choice of $u,v,w$,
 (\ref{ycw11}) can be written as
  \be\label{ycw15}
 \alpha=4\mu^{-1}(\mu/4+\epsilon\delta^2\mu^{-1}-\epsilon c^2)h_{\mu},\ \
     \beta=4\mu^{-\frac{3}{2}}\sqrt{\mu/4+\epsilon\delta^2\mu^{-1}-\epsilon c^2}\
     c_0,
  \ee
 and {\bf K} in Theorem \ref{th2}(iic) has the form
 \be\label{ycw16}
{\bf K}=\frac{6\mu(a_1^2-4\epsilon)(1-\epsilon
s^2)^2c^2+(\mu^2+4\epsilon \delta^2)(a_1\epsilon s^3+6\epsilon
s^2+3a_1s+2)(1+a_1s+\epsilon
s^2)}{128\mu^{-2}(\mu/4+\epsilon\delta^2\mu^{-1}-\epsilon
c^2)^3(1+a_1s+\epsilon s^2)^4}.
 \ee
 \een
\end{thm}

For a square metric $F=(\alpha+\beta)^2/\alpha$, we have
$a_1=2,\epsilon=1$, and then (\ref{ycw16}) becomes (\ref{ycw13})
and {\bf K} is bounded by (\ref{ycw14}). In other cases, {\bf K}
might not have positive lower bound.

The metric in Theorem \ref{th3} can be equivalently considered as
a subclass of $\phi(s)$ determined by (\ref{y1}) with
$k_1,k_2,k_3$ satisfying certain relation (see Remark \ref{rem72}
below).

\bigskip

It might be true that $\beta$ is closed for a regular
$(\alpha,\beta)$-metric (not of Randers type) of scalar flag
curvature in dimension $n\ge 3$. But it seems hard to be proved.
On the other hand, if an  $(\alpha,\beta)$-metric $F$ is singular
with $\phi(0)=0$ or $\phi(0)$ not defined, $\beta$ is generally
not closed even if $F$ is locally projectively flat or of constant
flag curvature (cf. \cite{SYa1} \cite{Y3} \cite{Y4}). We
conjecture that an $(\alpha,\beta)$-metric (regular or singular,
and not of Randers type or Kropina type) must be locally
projectively flat if it is of scalar flag curvature in dimension
$n\ge 3$.

\section{Preliminaries}\label{sec2}

  In local coordinates, the geodesics of a Finsler metric
  $F=F(x,y)$ are characterized by
   $$\frac{d^2 x^i}{d t^2}+2G^i(x,\frac{d x^i}{d t})=0,$$
where
 \beq \label{G1}
 G^i:=\frac{1}{4}g^{il}\big \{[F^2]_{x^ky^l}y^k-[F^2]_{x^l}\big \}.
 \eeq
For a Finsler metric $F$, the Riemann curvature $R_y=R^i_{\
k}(y)\frac{\pa}{\pa x^i}\otimes dx^k$ is defined by
 \be\label{y7}
 R^i_{\ k}:=2\frac{\pa G^i}{\pa x^k}-y^j\frac{\pa^2G^i}{\pa x^j\pa
 y^k}+2G^j\frac{\pa^2G^i}{\pa y^j\pa y^k}-\frac{\pa G^i}{\pa y^j}\frac{\pa G^j}{\pa
 y^k}.
 \ee
The Ricci curvature is the trace of the Riemann curvature, ${\bf
Ric}:=R^m_{\ m}$.

A Finsler metric $F$ is said to be projectively flat in $U$, if
 there is a local coordinate system $(U,x^i)$ such that
  $G^i=Py^i$,
  where $P=P(x,y)$ is called the  projective factor.

In projective geometry, the Weyl curvature and the Douglas
curvature play a very important role. We first introduce their
definitions. Put
 $$A^i_{\ k}:=R^i_{\ k}-R\delta^i_k, \ \ \  R:=\frac{R^m_{\
 m}}{n-1}.$$
Then the Weyl curvature $W^i_{\ k}$ are defined by
 \be\label{y8}
 W^i_{\ k}:=A^i_{\ k}-\frac{1}{n+1}\frac{\pa A^m_{\ k}}{\pa
 y^m}y^i.
 \ee
The Douglas curvature $D^{\ i}_{h \ jk}$ are defined by
 $$D^{\ i}_{h \ jk}:=\frac{\pa^3}{\pa y^h\pa y^j\pa
 y^k}\big(G^i-\frac{1}{n+1}G^m_my^i\big), \ \ \  G^m_m:=\frac{\pa G^m}{\pa y^m}.$$
The Weyl curvature and the Douglas curvature both are projectively
invariant. A Finsler metric is called a Douglas metric if $D^{\
i}_{h \ jk}=0$. A Finsler metric is of scalar flag curvature if
and only if $W^i_{\ k}=0$. It is known  that a Finsler metric with
the dimension $n\ge 3$ is locally projectively flat if and only if
$W^i_{\ k}=0$ and $D^{\ i}_{h \ jk}=0$ (\cite{Ma1}).

\bigskip

Now for the computation of $(\alpha,\beta)$-metrics, we first show
some definitions and conventions. For a Riemannian  $\alpha
=\sqrt{a_{ij}y^iy^j}$ and a $1$-form $\beta = b_i y^i $, let
 $$r_{ij}:=\frac{1}{2}(b_{i|j}+b_{j|i}),\ \ s_{ij}:=\frac{1}{2}(b_{i|j}-b_{j|i}),\ \
 r^i_{\ j}:=a^{ik}r_{kj},\ \  s^i_{\ j}:=a^{ik}s_{kj},$$
 $$q_{ij}:=r_{im}s^m_{\ j}, \ \ t_{ij}:=s_{im}s^m_{\ j},\ \
 r_j:=b^ir_{ij},\ \  s_j:=b^is_{ij},$$
 $$
  q_j:=b^iq_{ij}, \ \ r_j:=b^ir_{ij},\ \ t_j:=b^it_{ij},
 $$
 where we define $b^i:=a^{ij}b_j$, $(a^{ij})$ is the inverse of
 $(a_{ij})$, and $\nabla \beta = b_{i|j} y^i dx^j$  denotes the covariant
derivatives of $\beta$ with respect to $\alpha$. Here are some of
our conventions in the whole paper. For a general tensor $T_{ij}$
as an example, we  define $T_{i0}:=T_{ij}y^j$ and
$T_{00}:=T_{ij}y^iy^j$, etc. We use $a_{ij}$ to raise or lower the
indices of a tensor.

An $(\alpha,\beta)$-metric is a Finsler metrics defined by a
Riemann metric $\alpha=\sqrt{a_{ij}(x)y^iy^j}$ and a 1-form
$\beta=b_i(x)y^i$ as follows:
 $$F=\alpha \phi(s),\ \ s=\beta/\alpha,$$
 where $\phi(s)>0$ is a $C^{\infty}$ function on
$(-b_o,b_o)$. It is proved in \cite{Shen02} that an
$(\alpha,\beta)$-metric is regular if and only if
 \be\label{gjcw021}
 \phi(s)-s\phi'(s)+(b^2-s^2)\phi''(s)>0,\ \ (|s|\le b<b_o).
 \ee
 In this paper, an $(\alpha,\beta)$-metric is always assumed to be
 regular.

  Let $F=\alpha \phi(s),
s=\beta/\alpha$, be an $(\alpha,\beta)$-metric, and then  by
(\ref{G1}), the spray coefficients $G^i$ of $F$
 are given by
  \be\label{y005}
  G^i=G^i_{\alpha}+\alpha Q s^i_0+\alpha^{-1}\Theta (-2\alpha Q
  s_0+r_{00})y^i+\Psi (-2\alpha Q s_0+r_{00})b^i,
  \ee
where $G^i_{\alpha}$ denote the spray coefficents of $\alpha$ and
 $$
  Q:=\frac{\phi'}{\phi-s\phi'},\ \
  \Theta:=\frac{Q-sQ'}{2\Delta},\ \
  \Psi:=\frac{Q'}{2\Delta},\ \ \Delta:=1+sQ+(b^2-s^2)Q'.
 $$

For a general $(\alpha,\beta)$-metric $F=\alpha \phi(s),
s=\beta/\alpha$, we can use Maple programs, by aid of
(\ref{y005}), (\ref{y7}) and (\ref{y8}), to compute its Riemann
curvature $R^i_{\ k}$, Ricci curvature ${\bf Ric}$ and Weyl
curvature $W^i_{\ k}$. But these expressions are too complicated.
One can see the very long expressions of the Rimann curvature and
Ricci curvature in \cite{CST}. So for the briefness in this paper,
we won't write out the expression of the Weyl curvature for the
$(\alpha,\beta)$-metrics discussed in our paper.

We can show some details in the computation of the Weyl curvature
$W^i_{\ k}$ for an $(\alpha,\beta)$-metric $F=\alpha \phi(s),
s=\beta/\alpha$. The spray $G^i$ of $F$ can be written in the form
 $$G^i=\bar{G}^i+Py^i+Q^i=\widetilde{G}^i+Py^i,\ \ \ \widetilde{G}^i:=G^i_{\alpha}+Q^i,$$
where
 $$P:=\alpha^{-1}\Theta (-2\alpha Q
  s_0+r_{00}), \ \ Q^i:=\alpha Q s^i_0+\Psi (-2\alpha Q s_0+r_{00})b^i.$$
$\widetilde{G}^i$ define a spray $\widetilde{G}$ which is
projectively equivaent to the spray $G$ of $F$. The Weyl curvature
of $\tilde{G}$ is defined in the same way. $Q^i$ also define a
tensor $Q$.
 Since the
Weyl curvature is projectively invariant, we have $W^i_{\
k}=\widetilde{W}^i_k$, and then we can change to the computation
of $\widetilde{W}^i_k$ for $\widetilde{G}$. First we have
 $$\widetilde{R}^i_{\ k}=\bar{R}^i_{\ k}+H^i_{\ k}, \ \ \widetilde{R}=R+H,$$
 where
  \be\label{y10}
  H^i_{\
  k}:=2Q^i_{\ |k}-Q^i_{\ |j\cdot k}y^j+2Q^jQ^i_{\cdot j\cdot k}-Q^i_{\cdot j}Q^j_{\cdot k},\ \
  H:=\frac{H^m_{\ m}}{n-1}.
  \ee
where $Q^i_{\ |k}$  are the covariant derivatives of $Q$ with
respect to the spray $G_{\alpha}$ of $\alpha$. We have
 $$\widetilde{A}^i_{\ k}=\widetilde{R}^i_{\ k}-\widetilde{R}\delta^i_k=\bar{A}^i_{\ k}+H^i_{\ k}-H\delta^i_k.$$
Thus we get
 $$W^i_{\ k}=\widetilde{W}^i_k=\widetilde{A}^i_{\ k}-\frac{1}{n+1}\frac{\pa \widetilde{A}^m_{\ k}}{\pa
 y^m}y^i=\bar{W}^i_{\ k}+\Theta^i_{\ k},$$
 where
 \be\label{y11}
  \Theta^i_{\ k}:=(H^i_{\ k}-H\delta^i_k)-\frac{1}{n+1}\frac{\pa}{\pa y^m}\big(H^m_{\ k}-H\delta^m_k\big)y^i.
  \ee

As seen from above, we can use (\ref{y10}) and (\ref{y11}) to
compute the Weyl curvature $W^i_{\ k}$. In our final result, we
are mainly concerned about the computation of the following terms
 $$r_{ij}, \ \ r_{ij|k},\ \ q_{ij},\ \ t_{ij},\ \ s_{ij}.$$

We show some ideas in dealing with the equation $W^i_{\ k}=0$. We
first multiply $W^i_{\ k}=0$ by some suitable term and then
 it can be written in the following form
  \be\label{ycw}
  f_0(s)+f_1(s)\alpha+f_2(s)\alpha^2+\cdots +f_m(s)\alpha^m=0,
  \ee
  where $f_i(s)$'s are polynomials of $s$ with coefficients being
  homogenous polynomials in $(y^i)$.
 The key point is to
choose some suitable polynomials in $s$ to divide our equations,
and then get some answers by isolating rational and irrational
terms.

\section{Scalar Flag Curvature}
In this section, we study the properties of some class of
$(\alpha,\beta)$-metrics which are of scalar flag curvature. We
have the following two Theorems.

\begin{thm}\label{th1}
  Let $F=\alpha\phi(\beta/\alpha)$ be an
  $(\alpha,\beta)$-metric on an  $n(\ge 3)$-dimensional
 manifold $M$, where $\alpha =\sqrt{a_{ij}(x)y^iy^j}$ is Riemannian and $\beta = b_i(x) y^i$ is a
 1-form, and $\phi(s)$ satisfies (\ref{y1}) with  $\phi(0)=1$ and $k_2\ne k_1k_3$. Assume $\beta$ is closed. Then $F$
    is of scalar flag curvature     if and only if
 the Riemann  curvature $\bar{R}^i_{\ k}$ of $\alpha$ and the covariant
    derivatives $b_{i|j}$ of $\beta$ with respect to $\alpha$ satisfy the following equations
   \beq
    b_{i|j}&=&\tau
     \big\{(1+k_1b^2)a_{ij}+(k_2b^2+k_3)b_ib_j\big\}, \ \ \ (b^2:=||\beta||_{\alpha}^2),\label{y2}\\
   \bar{R}^i_{\ k}& =&\lambda
   (\alpha^2\delta^i_k-y^iy_k)+\eta\big(\beta^2\delta^i_k+\alpha^2b^ib_k-\beta
   b^iy_k-\beta b_ky^i),\label{y3}\\
  \tau_{x^i}&=&q b_i, \label{y4}
   \eeq
where $\lambda=\lambda(x),\tau=\tau(x)$ are scalar functions on
$M$ and $\eta,q$ are defined by
 \be\label{y5}
 \eta:=\big\{k_1^2+k_2-2k_1k_3-k_1(k_2-k_1^2)b^2\big\}\tau^2+k_1\lambda,\ \ \
 q:=(k_3-2k_1-k_1^2b^2)\tau^2-\lambda.
 \ee
In this case, $F$ is locally projectively flat, and the scalar
flag curvature ${\bf K}$ is given by
 \beq
32\phi^2{\bf
K}&=&f\phi'\phi^{-1}\Big\{\big[24f\phi^{-1}\phi'+(3-3f+hs^2)^2s^{-3}b^2-16hs\big]\tau^2+16\lambda
s\Big\}\nonumber\\
&&+\big[8(2ghs^2+12f-3g^2)-g(3+hs^2-3f)^2s^{-2}b^2\big]s^{-2}\tau^2-16\lambda
g,\label{y6}
 \eeq
 where $f,g,h$ are defined by
  \be\label{y007}
 f:=1+(k_1+k_3)s^2+k_2s^4, \ \ g:=k_2s^4-k_1s^2-2, \ \
 h:=3k_2s^2-k_1+3k_3.
  \ee
\end{thm}

 If $\tau=0$ in
 (\ref{y2}), we have $\beta=0$ and $\alpha$ is of constant sectional curvature, or $\beta\ne 0$ and $\alpha$ is
 flat ($\lambda=\eta=\tau=q=0$). Further, $\lambda$ and $\tau$ in (\ref{y2})--(\ref{y4}) can also be
 determined (see (\ref{gjcw63}), (\ref{gjcw66}) and (\ref{yy084}) below). The
 inverse of Theorem \ref{th1} is also true, which is shown as
 follows.

\begin{thm}\label{th31}
  Let $F=\alpha\phi(\beta/\alpha)$ be an
  $(\alpha,\beta)$-metric on an  $n(\ge 3)$-dimensional
 manifold $M$, where $\alpha =\sqrt{a_{ij}(x)y^iy^j}$ is Riemannian and $\beta = b_i(x) y^i$ is a
 1-form, and $\phi(0)=1$. Assume $F$ is not of Randers type, $\beta$ is not parallel with respect to
 $\alpha$.  Then $\phi(s)$ must
    satisfy (\ref{y1}), provided that $F$
    is of scalar flag curvature and (\ref{y2})--(\ref{y5}) hold with $k_2\ne k_1k_3$.
\end{thm}

\subsection{Proof of Theorem \ref{th1}}

We will deal with the equation $W^i_{\ k}=0$ step by step. By the
method described in  Section \ref{sec2}, we get a formula for
$W^i_{\ k}$ which is expressed in terms of  the covariant
derivatives
 of $\beta$ with respect to $\alpha$ and the Riemann curvature of $\alpha$.

Since $F$ in Theorem \ref{th1} is regular, we have the following
lemma.

\begin{lem}\label{lem1}
 Let $F=\alpha\phi(\beta/\alpha)$ be an
 $(\alpha,\beta)$-metric, where $\phi(s)$ is given by (\ref{y1}).
 Then we have
  $$1+k_1b^2>0,\ \ 1+(k_1+k_3)b^2+k_2b^4>0.$$
\end{lem}

Since $\beta$ is assumed to be closed, we have
$s_{ij}=0,q_{ij}=0,t_{ij}=0$. Plugging them into $W^i_{\ k}$ and
multiplying $W^i_{\ k}=0$ by
 $$4(n^2-1)\big[\phi-s\phi'+(b^2-s^2)\phi''\big]^5\alpha^4$$
and we get an equation  denoted by $Eq_0=0$. By  (\ref{y1}) we can
get $\phi^{(i)}$ $(2\le i\le 5)$ expressed by $\phi,\phi'$. Plug
them into $Eq_0=0$ and then multiply $Eq_0=0$ by
 $$\big[1+(k_1+k_3)s^2+k_2s^4\big]^5.$$
By this way, we have
 $$
 4(\phi-s\phi')^5Eq_1=0.
 $$
It is surprising that $Eq_1$ is independent of $\phi$ and by
$Eq_1=0$ we have
 \beq\label{gjcw24}
 0&=&24(n-2)(k_2-k_1k_3)^3(\alpha^2
 b_k-\beta y_k)y^i\beta^3\big[\alpha^4+(k_1+k_3)\alpha^2\beta^2+k_2\beta^4\big]^2r_{00}^2\nonumber\\
 &&+C^i_k\big[(1+k_1b^2)\alpha^2+(k_2b^2+k_3)\beta^2\big],
 \eeq
where $C^i_k$ are homogeneous polynomials in $(y^i)$.

 \begin{lem}\label{lem2}
$(1+k_1b^2)\alpha^2+(k_2b^2+k_3)\beta^2$ cannot be divided by
$\alpha^2 b_k-\beta y_k$ for all $k$.
 \end{lem}

 {\it Proof :}  Otherwise, for some scalar functions $f_k=f_k(x)$ we have
 $$\alpha^2 b_k-\beta
 y_k=f_k\big[(1+k_1b^2)\alpha^2+(k_2b^2+k_3)\beta^2\big].$$
 Then we have
  $$ b^2\alpha^2-\beta^2=f\big[(1+k_1b^2)\alpha^2+(k_2b^2+k_3)\beta^2\big], \ \ f:=b^kf_k,$$
which implies
 $1+k_1b^2+(k_2b^2+k_3)b^2=0.$
 This is a contradiction by Lemma \ref{lem1}.   \qed

\begin{lem}\label{lem3}
$(1+k_1b^2)\alpha^2+(k_2b^2+k_3)\beta^2$ cannot be divided by
$\alpha^4+(k_1+k_3)\alpha^2\beta^2+k_2\beta^4$, provided that
$k_2\ne k_1k_3$.
 \end{lem}

{\it Proof :} We can prove it in two cases: $k_2b^2+k_3=0$ and
$k_2b^2+k_3\ne 0$. We need Lemma \ref{lem1} and the fact $k_2\ne
k_1k_3$.   \qed

\bigskip

By (\ref{gjcw24}), Lemma \ref{lem2} and Lemma \ref{lem3} we have
(\ref{y2}) for some scalar function $\tau=\tau(x)$.  Then by
(\ref{y2}) we can obtain the expressions of the following
quantities:
 $$
r_{00},\ r_i,\ r^m_m,\ r, \ r_{00|0}, \ r_{0|0}, \ r_{00|k}, \
r_{k0|m}, \ r_{k|0}, \ etc.
 $$
 For example, we have
 $$r_{0|0}=\big[1+(k_1+k_3)b^2+k_2b^4\big]\Big\{\big[(1+k_1b^2)\alpha^2+(2k_1+3k_3+5k_2b^2)\beta^2\big]\tau^2+\tau_0\beta\Big\}.$$
Now plug all the above quantities into (\ref{gjcw24}) and then we
can get the Weyl curvature $\bar{W}_{ik}:=a_{im}\bar{W}^m_{\ k}$
of $\alpha$. We will discuss it under two cases.

\bigskip

\noindent{\bf Case I:} Assume $k_1\ne 0$. We have
 \beq\label{gjcw25}
 \bar{W}_{ik}&=&\frac{k_1}{n-1}b^m\omega_m(\alpha^2a_{ik}-y_iy_k)-\frac{(k_1+k_2b^2)\omega_0-k_2b^m\omega_m\beta}{n-1}\ \beta
  a_{ik}\nonumber\\
 &&+\frac{1}{n-1}\Big\{\frac{k_1+k_2b^2}{n+1}\big[(2n-1)\beta \omega_k-(n-2)\omega_0b_k\big]
 -k_2b^m\omega_m\beta b_k\Big\}y_i\nonumber\\
 &&+\omega_0b_i(k_1y_k+k_2\beta b_k)-(k_1\alpha^2+k_2\beta^2)b_i\omega_k,
 \eeq
where $\tau_i:=\tau_{x^i}$ and
 \be\label{gjcw025}
 \omega_i:=\tau_i+\Big(k_1+k_3+k_2b^2-\frac{k_2}{k_1}\Big)\tau^2b_i.
 \ee

\begin{lem}\label{lem4}
 (\ref{gjcw25}) $\Longleftrightarrow $ (\ref{y3}) and (\ref{y4}),
 where $q$ is defined by
  \be\label{gjcw0025}
 q:=-\frac{\eta}{k_1}-\big(k_1+k_3+k_2b^2-\frac{k_2}{k_1}\big)\tau^2.
  \ee
\end{lem}

{\it Proof} :  {\bf $\Longrightarrow$ :}
 By the definition of the Weyl curvature
$\bar{W}_{ik}$ of $\alpha$ we have
 \be\label{y52}
 \bar{W}_{ik}=\bar{R}_{ik}-\frac{1}{n-1}\bar{R}ic_{00}a_{ik}+\frac{1}{n-1}\bar{R}ic_{k0}y_i,
 \ee
where $\bar{R}_{ik}:=a_{im}\bar{R}^m_{\ k}$ and $\bar{R}ic_{ik}$
denote the Ricci tensor of $\alpha$. Using the fact
$\bar{R}_{ik}=\bar{R}_{ki}$ we get from  (\ref{y52})
 \be\label{y53}
\bar{W}_{ik}-\bar{W}_{ki}=\frac{1}{n-1}\big(\bar{R}ic_{k0}y_i-\bar{R}ic_{i0}y_k\big).
 \ee
  By (\ref{gjcw25}) we can get another expression of
  $\bar{W}_{ik}-\bar{W}_{ki}$. Thus by (\ref{gjcw25}) and
  (\ref{y53}) we have

  \be\label{y54}
 T_iy_k-T_ky_i+(n^2-1)(k_1\alpha^2+k_2\beta^2)(\omega_ib_k-\omega_kb_i)=0,
  \ee
 where we define
  \beqn
  T_i:&=&(n+1)\bar{R}ic_{i0}-(2n-1)(k_1+k_2b^2)\beta\omega_i\\
  &&+\big\{[(n^2+n-3)k_1+(n-2)k_2b^2]\omega_0+(n+1)k_2b^m\omega_m\beta\big\}b_i.
  \eeqn
Contracting (\ref{y54}) by $y^k$ we get
 \be\label{y55}
 \big[T_i+(n^2-1)k_1(\omega_i\beta-\omega_0b_i)\big]\alpha^2-T_0y_i+(n^2-1)k_2\beta^2(\omega_i\beta-\omega_0b_i)=0.
 \ee
Contracting (\ref{y55}) by $b^i$ we obtain
 \be\label{gjcw30}
\big[b^mT_m+(n^2-1)k_1(b^m\omega_m\beta-b^2\omega_0)\big]\alpha^2
+\big[(n^2-1)k_2\beta(b^m\omega_m\beta-b^2\omega_0)-T_0\big]\beta=0.
 \ee
 So by (\ref{gjcw30}) there is some scalar
function $\bar{\eta}=\bar{\eta}(x)$ such that
 \be\label{y56}
 T_0=(n^2-1)k_2\beta(b^m\omega_m\beta-b^2\omega_0)+(n+1)\bar{\eta}\alpha^2.
\ee
 Then by the definition of $T_i$ and (\ref{y56}) we have
  \beq\label{y57}
\bar{R}ic_{00}&=&\bar{\eta}\alpha^2-(n-2)\big[(k_1+k_2b^2)\omega_0-k_2b^m\omega_m\beta\big]\beta,\label{y57}\\
\bar{R}ic_{i0}&=&\bar{\eta}
y_i-(n-2)\big\{\frac{k_1+k_2b^2}{2}(\beta\omega_i+b_i
\omega_0)-k_2b^m\omega_m\beta b_i\big\}.\label{y057}
  \eeq
Now plugging (\ref{y57}) and (\ref{y057}) into (\ref{y55}) yields
\be\label{gjcw34}
 2(n+1)k_2A_i\beta+B_i\alpha^2=0,
\ee
 where $A_i$ and $B_i$ are defined by
 \beqn
 A_i:&=&(b^2\omega_0-b^m\omega_m\beta)y_i+\beta^2\omega_i-\beta\omega_0b_i,\\
 B_i:&=&k_2\big\{\big[2(n+1)b^m\omega_m\beta-(n-2)b^2\omega_0\big]b_i-(n+4)b^2\beta\omega_i\big\}+(n-2)k_1(\beta\omega_i-\omega_0b_i).
\eeqn
 If $k_2=0$, then by (\ref{gjcw34}) we have
  \be\label{gjcw35}
 \beta\omega_i-\omega_0b_i=0.
  \ee
If $k_2\ne0$, then by (\ref{gjcw34}) we have
 \be\label{gjcw36}
 A_i=f_i\alpha^2,
 \ee
where $f_i=f_i(x)$ are scalar functions. Contracting
(\ref{gjcw36}) by $y^i$ we get
 \be\label{gjcw37}
 f_i=b^2\omega_i-b^m\omega_mb_i
 \ee
Plugging (\ref{gjcw36}) and (\ref{gjcw37}) into (\ref{gjcw34})
gives
 \be\label{gjcw38}
 (n-2)(k_1+k_2b^2)(\beta\omega_i-\omega_0b_i)\alpha^2=0.
 \ee
 If $\tau= 0$, we can naturally find $\lambda,\eta,q$ such that
 (\ref{y2})--(\ref{y4}) hold.
So we may assume $\tau\ne 0$, and then by  (\ref{y1}), we have $b^2\ne constant$. So by (\ref{gjcw38}) we also
get (\ref{gjcw35}). Thus it follows from (\ref{gjcw35}) that
 \be\label{gjcw39}
 \omega_i=eb_i,
 \ee
for some scalar function $e=e(x)$. Now plugging (\ref{y57}),
(\ref{y057}) and (\ref{gjcw39}) into (\ref{gjcw25}) and
(\ref{y52}) we obtain
 \be\label{y59}
 \bar{R}_{ik} =\lambda
   (\alpha^2a_{ik}-y^iy_k)+\eta\big(\beta^2a_{ik}+\alpha^2b_ib_k-\beta
   b_iy_k-\beta b_ky_i),
 \ee
where we define
 \be\label{y60}
 \lambda:=\frac{k_1eb^2+\bar{\eta}}{n-1},\ \ \ \eta:=-k_1e.
 \ee
Clearly, (\ref{y59}) is just (\ref{y3}). It follows from
(\ref{gjcw025}), (\ref{gjcw39}) and (\ref{y60}) that
 \be
 \tau_i=qb_i=\Big\{-\frac{\eta}{k_1}-\big(k_1+k_3+k_2b^2-\frac{k_2}{k_1}\big)\tau^2\Big\}b_i,
 \ee
which implies (\ref{y4}) with $q$ given by (\ref{gjcw0025}).

{\bf $\Longleftarrow$ :} We verify that both sides of
(\ref{gjcw25}) are equal. By (\ref{y3}) we have (\ref{y57}) and
(\ref{y057}). Since (\ref{y52}) naturally holds, we plug
(\ref{y57}), (\ref{y057}) and (\ref{y3}) into (\ref{y52}) and then
we obtain the left side of (\ref{gjcw25}). By (\ref{y4}) and
(\ref{gjcw0025}) we get
 \be\label{gjcw45}
 \omega_i=eb_i=-\frac{\eta}{k_1}
 \ee
 from (\ref{gjcw025}).
 Then plugging (\ref{gjcw45}) into the right
side of (\ref{gjcw25}) we obtain the result equal to the left side
of (\ref{gjcw25}).   \qed

\bigskip

Now we show that the metric $F$ in Theorem \ref{th1} is locally
projectively flat if $F$ is of scalar flag curvature with $\beta$
being closed and the dimension $n\ge 3$. First, we have the
vanishing Weyl curvature $W^i_{\ k}=0$ since $F$ is of scalar flag
curvature. Second, by (\ref{y2}) $F$ is a Douglas metric
(\cite{LSS}). Therefore it follows from the result in \cite{Ma1}
that $F$ is locally projectively flat.

In the final we compute the scalar flag curvature and prove that
$\eta,q$ are given by (\ref{y5}). As shown above, we have
(\ref{y2})--(\ref{y4}) since $F$ is of scalar flag curvature. Plug
(\ref{y1})  and (\ref{y2})--(\ref{y4}) into the Riemann curvature
$R^i_{\ k}$ of $F$, where $q$ is defined by (\ref{gjcw0025}), and
then a direct computation gives
 \beq\label{gjcw46}
 R^i_{\ k} &=&K
 F^2(\delta^i_k-F^{-1}y^iF_{y^k})+\frac{\phi'(sF_{y^k}-\phi
 b_k)}{k_1\phi^2(\phi-s\phi')}\times \nonumber\\
 &&\big[\eta-\big\{k_1^2+k_2-2k_1k_3-k_1(k_2-k_1^2)b^2\big\}\tau^2-k_1\lambda\big]Fy^i,
 \eeq
 where the expression of $K=K(x,y)$ is omitted.
Since $F$ is of scalar flag curvature and $n\ge 3$, by
(\ref{gjcw46}) we must have
 \be\label{gjcw47}
\eta-\big\{k_1^2+k_2-2k_1k_3-k_1(k_2-k_1^2)b^2\big\}\tau^2-k_1\lambda=0.
 \ee
Thus we get $\eta$ given by (\ref{y5}). Plug $\eta$ given by
(\ref{y5}) into $K$, and then we obtain the scalar flag curvature
${\bf K}=K$ given by (\ref{y6}). By (\ref{gjcw0025}) and $\eta$ in
(\ref{y5}), we get $q$ given by (\ref{y5}).

\bigskip

\noindent{\bf Case II:} Assume $k_1= 0$. Then $k_2\ne 0$. We have
 \beq\label{gjcw48}
 \frac{\bar{W}_{ik}}{k_2}&=&\frac{(\tau^2+b^m\tau_m)\beta^2-b^2(\beta\tau_0+\tau^2\alpha^2)}{n-1}
 \delta_{ik}+(\tau_0b_k-\beta\tau_k)\beta b_i+\tau^2(\alpha^2 b_k-\beta y_k)b_i\nonumber\\
 &&+\frac{\big[(2n-1)\beta\tau_k-(n-2)\tau_0b_k\big]b^2+(n+1)\big[b^2\tau^2y_k-(\tau^2+b^m\tau_m)\beta
 b_k\big]}{n^2-1}\ y_i.
 \eeq

\begin{lem}
 (\ref{gjcw48}) $\Longleftrightarrow $ (\ref{y3}) and (\ref{y4}),
 where $q=q(x)$ is some scalar function and $\eta=k_2\tau^2$.
\end{lem}

{\it Proof} : We only show the final results of the necessity.
Assume (\ref{gjcw48}). By a similar steps as in Lemma \ref{lem4},
we have
 \beq
 \tau_i&=&qb_i,\label{gjcw49}\\
 \bar{R}_{ik} &=&\lambda
   (\alpha^2a_{ik}-y^iy_k)+k_2\tau^2\big(\beta^2a_{ik}+\alpha^2b_ib_k-\beta
   b_iy_k-\beta b_ky_i),\label{gjcw50}
 \eeq
where $q=q(x),\lambda=\lambda(x)$ are scalar functions.  By
(\ref{gjcw50}), $\eta$ in (\ref{y3}) is given by $\eta=k_2\tau^2$,
which is equal to the $\eta$ in (\ref{y5}) with $k_1=0$.   \qed

\bigskip

Now we compute the Riemann curvature $R^i_{\ k}$ of $F$ in this
case and prove that $\eta$ are given by (\ref{y5}) with $k_1=0$.
By (\ref{y1})  and (\ref{y2})--(\ref{y4}) with $\eta=k_2\tau^2$, a
direct computation gives
 \be\label{gjcw51}
 R^i_{\ k} =K
 F^2(\delta^i_k-F^{-1}y^iF_{y^k})-\frac{\phi'(sF_{y^k}-\phi
 b_k)}{\phi^2(\phi-s\phi')}
 (\lambda+q-k_3\tau^2)Fy^i,
 \ee
 where the expression of $K=K(x,y)$ is omitted.
Since $F$ is of scalar flag curvature and $n\ge 3$, by
(\ref{gjcw51}) we must have
 \be\label{gjcw52}
\lambda+q-k_3\tau^2=0.
 \ee
Thus we get $q$ given by (\ref{y5}) with $k_1=0$.

\subsection{Proof of Theorem \ref{th31}}

We will show that $\phi(s)$ must satisfy (\ref{y1}) by $W^i_{\
k}=0$ if (\ref{y2})--(\ref{y5}) hold with $k_2\ne k_1k_3$,
$\tau\ne 0$ and $\phi(s)\ne a_1s+\sqrt{1+ks^2}$ for any constants
$a_1,k$.

To determine $\phi$ by $W^i_{\ k}=0$, we need to choose a special
coordinate system at a point on $M$ as that in \cite{Shen02}. At a
fixed point $x_o$, make a change of coordinates: $(s,y^A)\mapsto
(y^1,y^A)$ by
 $$y^1=\frac{s}{\sqrt{b^2-s^2}}\bar{\alpha},\ \ y^A=y^A,$$
where $\bar{\alpha}=\sqrt{\sum_{A=2}^n (y^A)^2}$. Then
$$\alpha=\frac{b}{\sqrt{b^2-s^2}}\bar{\alpha},\ \
\beta=\frac{bs}{\sqrt{b^2-s^2}}\bar{\alpha}.$$

Using (\ref{y52}), plug (\ref{y2})--(\ref{y5}) into the Weyl
curvature $W^i_{\ k}$ of $F$, and then multiplying them by
 $$
 4(n^2-1)\alpha^4\big[\phi-s\phi'+(b^2-s^2)\phi''\big],
 $$
we get a family of ODEs about $\phi$. Here we only consider
$W^A_{\ B}=0$. Under the above coordinate system $(s,y^A)$,
$W^A_{\ B}=0$ can be written in the form
 $$
 (n+1)b^2\bar{\alpha}^2\big[\phi-s\phi'+(b^2-s^2)\phi''\big]Y_1\delta_{AB}+Y_2y_Ay_B=0,
 $$
where $Y_1,Y_2$ are two ODEs about $\phi$ independent of $(y_A)$.
Since $n\ge 3$, we have $Y_1=0,Y_2=0$.
  We only need to consider $Y_1=0$, which is a regular ODE of
order four with the leading term
 $$
 2\tau^2\big[1+k_1b^2+(k_3+k_2b^2)s^2\big]^3(\phi-s\phi')^2\big[\phi-s\phi'+(b^2-s^2)\phi''\big]\phi^{(4)}.
 $$ Plug
the Taylor expansion
 $$\phi=1+\sum_{i=1}^6 a_is^i+o(s^6)$$
 into $Y_1=0$, and let $p_i$ be the coefficients of $s^i$ in
 $Y_1$. By $p_0=0$ we have
  \be\label{gjcw53}
 f\tau^2+4\lambda(1+2a_2b^2)^3(2a_2-k_1)=0,
  \ee
where we omit the expression of $f=f(x)$.

\bigskip

\noindent{\bf Case I:} Assume $a_2=k_1/2$. By (\ref{gjcw53}) we
have
 \be\label{gjcw54}
 (3k_1^3+2k_1^2k_3-2k_1k_2+24a_4k_1-54a_3^2)b^2+3k_1^2-2k_2+2k_1k_3+24a_4=0.
 \ee
By (\ref{y1}) with $\tau\ne 0$ we have $b^2\ne constant$. So by
(\ref{gjcw54}), we easily get $a_3=0$. Since the $\phi$ determined
by (\ref{y1}) satisfies
 $$\phi(0)=1,\ \ \phi''(0)=k_1,\ \ \phi'''(0)=0,$$
  (\ref{y1}) is  the unique solution to $Y_1=0$ by the uniqueness
  of solution for regular ODEs.

\bigskip

\noindent{\bf Case II:} Assume $a_2\ne k_1/2$. Solving $\lambda$ from (\ref{gjcw53})  and
plugging it into $Y_1=0$ we get a regular ODE denoted by $Y_3=0$. We will show that this case will not occur.

 Solving $\lambda$ from (\ref{gjcw53})  and plugging it into
 $p_1=0$, using $b\ne constant$ we obtain
  \be\label{gjcw55}
  a_5=\frac{a_3(62a_2^2-36k_1a_2+8k_3a_2-12a_4-5k_1k_3+k_1^2+k_2)}{20(k_1-2a_2)},
  \ee
 and $a_3=0$ or
  \be \label{gjcw56}
  a_4=\frac{9a_3^2}{44(2a_2-k_1)}-\frac{a_2^2}{2}+\frac{a_2(k_3-k_1)}{66}+\frac{k_1(k_1-7k_3)}{132}+\frac{k_2}{22}.
  \ee

 {\bf Case II(1):} Assume $a_3=0$.  Solving $\lambda$ from (\ref{gjcw53})  and plugging it
 into $p_2=0$, using $b\ne constant$ and $a_2\ne k_1/2$ we obtain from $a_3=0$ and (\ref{gjcw55})
  \be\label{gjcw57}
  a_4=-\frac{a_2^2}{2},\ \  a_6=-\frac{a_2^2k_3+3a_2^3+11a_2a_4+2a_4k_3}{5}.
  \ee
  In this case, we can verify that
   $$\phi(s)=a_1s+\sqrt{1+2a_2s^2}$$
   is the unique solution of $Y_3=0$.  This contradicts with our
   assumption on $\phi$.

 {\bf Case II(2):} Assume $a_3\ne 0$. Solving $\lambda$ from (\ref{gjcw53})  and plugging it into
 $p_1=0$, by (\ref{gjcw55}) and (\ref{gjcw56}) we obtain
  $$
  a_3^2=\frac{(k_1-2a_2)(2a_2k_3-2a_2k_1+k_1^2+8k_2-9k_1k_3)}{189},
  $$
  and $k_1=k_3$ or
   $$a_2=\frac{k_1^2-11k_1k_3+10k_2}{2(k_3-k_1)}.$$
  In both cases, we will immediately show that $k_2=k_1k_3$, which
  is a contradiction with our assumption.          \qed

\section{Deformations on $(\alpha,\beta)$-metrics}\label{sec4}

 In Theorem \ref{th1}, (\ref{y2})--(\ref{y5}) are necessary and
 sufficient conditions for an $n(\ge 3)$-dimensional $(\alpha,\beta)$-metric which is not
 of Randers type satisfying (\ref{y1}) to be locally projectively flat. In \cite{Shen1},
 Shen gives another characterization, and then in \cite{Yu} \cite{Yu1}, Yu finds a
 deformation to obtain the local structure based on Shen's result.
 In this section, we will give the local structure using
 (\ref{y2})--(\ref{y5}) by a similar deformation.

 Let $u=u(t),v=v(t),w=w(t)$ satisfy the following ODEs:
 \beq
u'&=&\frac{v-k_1u}{1+(k_1+k_3)t+k_2t^2},\label{gjcw58}\\
v'&=&\frac{u(k_2u-k_3v-2k_1v)+2v^2}{u[1+(k_1+k_3)t+k_2t^2]},\label{gjcw59}\\
w'&=&\frac{w(3v-k_3u-2k_1u)}{2u\big[1+(k_1+k_3)t+k_2t^2\big]}.\label{gjcw60}
 \eeq
 Let $\alpha$ and $\beta$ satisfy (\ref{y2})--(\ref{y5}), and
define a new Riemann metric $h=\sqrt{h_{ij}(x)y^iy^j}$ and a new
1-form $\rho=p_i(x)y^i$ by
 \be\label{gjcw61}
 h:=\sqrt{u\alpha^2+v\beta^2},\ \ \rho:=w\beta,
 \ee
where $u=u(b^2)\ne 0,v=v(b^2),w=w(b^2)\ne 0$ which are determined
by (\ref{gjcw58})--(\ref{gjcw60}). We will show in the following
that $h$ is of constant sectional curvature and $\rho$ is a closed
1-form which is conformal with respect to $h$.

Put
 $$p:=||\rho||_{h}.$$
We have
 \be\label{gjcw62}
 p^2=\frac{w^2b^2}{u+vb^2}.
 \ee
By (\ref{y2}), a direct computation shows that the sprays $G^i_h$
of $h$ and $G^i_{\alpha}$ of $\alpha$ satisfy
 \be\label{gjcw64}
G^i_h=G^i_{\alpha}+\tau\Big\{\frac{1}{2}(k_1\alpha^2+k_2\beta^2)
b^i-\frac{(k_1u-v)\beta}{u}y^i\Big\},
 \ee

 By (\ref{y2}) and (\ref{gjcw64}) we can directly get
  \be\label{gjcw63}
 p_{i|j}=\frac{w\tau}{u}h_{ij} \ (=-2ch_{ij}),
  \ee
  where the covariant derivatives are taken with respect to $h$, and the
  scalar function $c=c(x)$ can be determined (see (\ref{yy084}) below).  Now (\ref{gjcw63})
 implies that $\rho$ is a closed conformal 1-form with respect to $h$.

By (\ref{gjcw64}) and using (\ref{y2})--(\ref{y5}), we obtain
 \be\label{gjcw65}
 \widetilde{R}^i_{\ k}=\frac{\lambda
u+(k_1^2ub^2+2k_1u-v)\tau^2}{u^2}(h^2\delta^i_k-y^i\widetilde{y}_k),
 \ee
where $\widetilde{R}^i_{\ k}$ are the Riemann curvatures of $h$
 and $\widetilde{y}_k:=h_{km}y^m$. It follows from (\ref{gjcw65})
that $h$ is of constant sectional curvature. We put it as $\mu$,
and then we obtain
 \be\label{gjcw66}
 \lambda=\mu u-\frac{k_1u(2+k_1b^2)-v}{u}\ \tau^2.
\ee
 In some local coordinate system we may put $h=h_{\mu}$ in the
form
 \be\label{gjcw67}
 h_{\mu}=\frac{\sqrt{(1+\mu |x|^2)|y|^2-\mu\langle
 x,y\rangle^2}}{1+\mu|x|^2},
 \ee
and then by (\ref{gjcw63}) and (\ref{gjcw67}) we obtain the 1-form
$\rho=p_iy^i$ given by
 \be\label{gjcw68}
 p_i=\frac{(k-\mu\langle a,x\rangle)
 x^i+(1+\mu|x|^2)a^i}{(1+\mu|x|^2)^{\frac{3}{2}}},\ \ \
 p^i=\sqrt{1+\mu|x|^2}(kx^i+a^i).
 \ee
where $k$ is a constant and $a=(a^i)$ is a constant vector, and
$p_i=h_{im}p^m$.  By (\ref{gjcw68}) we have
 \be\label{gjcw69}
 p^2=||\rho||^2_h=|a|^2+\frac{k^2|x|^2+2k \langle a,x\rangle -\mu \langle
 a,x\rangle^2}{1+\mu |x|^2}.
\ee

If we choose a triple $(u,v,w)$ determined by
(\ref{gjcw58})--(\ref{gjcw60}), then by the above discussion, we
can obtain the local expressions of $\alpha$ and $\beta$ by
(\ref{gjcw61}).

\begin{rem}
We can have different suitable choices of $u,v,w$ satisfying
(\ref{gjcw58})--(\ref{gjcw60}). For a square metric
$F=(\alpha+\beta)^2/\alpha$, the triple $(u,v,w)$ can be
 chosen as (\cite{SYa})
  $$u=(1-b^2)^2, \ \ v=0, \ \ w=\sqrt{1-b^2}.$$
  For the general case in Theorem \ref{th1}, we may choose the
  triple $(u,v,w)$ as (\cite{Yu} \cite{Yu1})
   \be\label{gjcw071}
   u=e^{2\sigma}, \ \ v=(k_1+k_3+k_2b^2)u, \ \ w=\sqrt{1+(k_1+k_3)b^2+k_2b^4}\ e^{\sigma},
   \ee
   where $\sigma$ is defined by
 $$2\sigma:=\int_0^{b^2}\frac{k_2t+k_3}{1+(k_1+k_3)t+k_2t^2}dt.$$
\end{rem}

\section{Constant flag curvature}\label{sec5}

In this section, we consider the classification in Theorem
\ref{th1} when $F$ is of constant flag curvature. The following
corollary has been proved in \cite{LS} \cite{Y2} in a different
way.

\begin{cor}\label{cor51}
 Let $F=\alpha\phi(\beta/\alpha)$ be an
  $(\alpha,\beta)$-metric on an  $n(\ge 2)$-dimensional
 manifold $M$, where $\phi(s)$ satisfies (\ref{y1}) with
 $\phi(0)=1$ and $k_2\ne k_1k_3$. Suppose $F$ is locally projectively flat with constant flag
 curvature $K$ and $\beta$ is not parallel with respect to $\alpha$.
 Then $F$ must be in the following form
 \be\label{F}
F=\frac{(\sqrt{\alpha^2+k\beta^2}+\epsilon\beta)^2}{\sqrt{\alpha^2+k\beta^2}},
\ee
 where $k$ and $\epsilon\ne 0$ are
     constant. In this case, we have $K=0$ and
     $k=k_1-\phi'(0)^2/2$.
\end{cor}

{\it Proof :} We can use Theorem \ref{th1} to prove it in case of
$n\ge 3$, and for the case $n=2$, see \cite{Y2}. Since $F$ is of
constant flag curvature, {\bf K} given by (\ref{y6}) is a
constant. Put
 $$\phi(s)=1+a_1s+\sum_{i=2}^{\infty}a_is^i.$$
Then by (\ref{y1}),  we can express all $a_i$'s ($i\ge 2$) in
terms of $k_1,k_2,k_3$. Multiply (\ref{y6}) by $\phi^2$ and we get
an equation. Let $p_i$ be the coefficients of $s^i$. Firstly, by
$p_0=0$, we get
 \be\label{gjcw72}
 K=\lambda+\big(k_1^2b^2+k_1+\frac{3}{4}a_1^2\big)\tau^2.
 \ee
We show $a_1\ne 0$. If $a_1=0$, then plugging $a_1=0$ and
(\ref{gjcw72}) into $p_2=0$ yields
 $$12\tau^2(k_2-k_1k_3)=0.$$
Since $k_2\ne k_1k_3$, we get $\tau=0$ on the whole $M$. Thus by
(\ref{y2}), $\beta$ is parallel with respect to $\alpha$. So
$a_1\ne 0$. Now substitute (\ref{gjcw72}) into $p_1=0$ and then
using $a_1\ne 0$ we obtain
 \be\label{gjcw73}
 \lambda=-(k_1^2b^2+k_3+2a_1^2)\tau^2.
 \ee
Next plugging (\ref{gjcw72}) and (\ref{gjcw73}) into $p_3=0$ and
using $a_1\ne 0$ and $\tau\ne 0$ we get
 \be\label{gjcw74}
 k_2=-a_1^4+\frac{3}{5}(k_1-k_3)a_1^2+\frac{1}{5}(k_1k_3+2k_1^2+2k_3^2).
 \ee
Then similarly, by (\ref{gjcw72})--(\ref{gjcw74}) and $p_4=0$ we
have
 $$k_3=k_1-a_1^2,\ k_1-\frac{5}{4}a_1^2,\
 -k_1+\frac{10}{3}a_1^2.$$
 If
 $$k_3=-k_1+\frac{10}{3}a_1^2,$$
 then plugging (\ref{gjcw72})--(\ref{gjcw74}) and $k_3$ into
 $p_5=0$ yields
  $$k_1=\frac{5}{12}a_1^2,\ \frac{13}{6}a_1^2,\
  \frac{55}{24}a_1^2.$$
It can be easily verify that if
 $$k_3=k_1-a_1^2, \ \ {\rm or} \ \ k_3=-k_1+\frac{10}{3}a_1^2\ \ {\rm and} \  \ k_1=\frac{5}{12}a_1^2,\
 \frac{13}{6}a_1^2,$$
  then we have $k_2=k_1k_3$. Therefore, we have
   \be\label{gjcw75}
 k_3=k_1-\frac{5}{4}a_1^2, \ \ {\rm or} \ \  k_3= -k_1+\frac{10}{3}a_1^2\ \ {\rm and} \ \ k_1=\frac{55}{24}a_1^2.
 \ee
The second case in (\ref{gjcw75}) is a special case of the first
case in (\ref{gjcw75}). So by (\ref{gjcw74}) and (\ref{gjcw75}),
we have
 \be\label{gjcw76}
 k_2=\frac{3}{8}a_1^4-\frac{5}{4}k_1a_1^2+k_1^2, \ \
 k_3=k_1-\frac{5}{4}a_1^2.
 \ee

Now by (\ref{gjcw72}), (\ref{gjcw73}) and (\ref{gjcw76}) we get
$K=0$. Plug (\ref{gjcw76}) into (\ref{y1}) and solving the ODE we
obtain (\ref{F}) with
 $$k:=k_1-a_1^2/2, \ \ \epsilon:=a_1/2.$$

\section{Proof of Theorem \ref{th2}}

To prove Theorem \ref{th2}, we need the following lemma.

\begin{lem} (\cite{SYa})\label{lem61}
 Let $\alpha=\sqrt{a_{ij}y^iy^j}$ be an $n$-dimensional
Riemannian metric of constant sectional curvature $\mu$ and
$\beta=b_iy^i$ is a 1-form on $M$. If $\beta$ satisfies
 \be\label{yy079}
 r_{ij}=-2ca_{ij},
 \ee
 where $c=c(x)$ is a scalar function on $M$, then the quantity
 $$
 |\nabla c|^2_{\alpha}+\mu c^2
 $$
  is a constant in case of $n\ge 3$, where $\nabla c$ is
 the gradient of $c$ with respect to $\alpha$. If $n\ge 2$ and $M$ is compact without boundary,
  then $c=0$ if $\mu<0$ and $c=constant$ if
 $\mu=0$.
 \end{lem}

\bigskip

{\it Proof of Theorem \ref{th2}} :

\bigskip

\noindent {\bf Case I:} In case of $n\ge 3$, suppose $F$ is of
constant flag
    curvature and $\beta$ is closed. Then $F$ is locally projectively
flat with constant flag curvature by Theorem \ref{th01}.

If $\tau=0$ in (\ref{y2}) on the whole $M$, then $\beta$ is
parallel and $\alpha$ is of constant sectional curvature on $M$.
So if $\alpha$ is flat, then $F$ is flat-parallel; if $\alpha$ is
of non-zero sectional curvature, then $\beta=0$ on $M$ and so
$F=\alpha$ is Riemannian.

If $\tau\ne 0$ at a point on  $M$, then by the proof in Corollary
\ref{cor51}, we have $a_1:=\phi'(0)\ne 0$, and
 $$
 \phi(s)=\frac{(\sqrt{1+ks^2}+\epsilon s)^2}{\sqrt{1+ks^2}},
 $$
and then $F$ is given by (\ref{F}). If (\ref{y0010}) holds, then
$\sqrt{\alpha^2+k\beta^2}$ is a (positive definite) Riemann
metric. So $F$ in (\ref{F}) is a square metric which can be
written as
 $$F=\frac{(\widetilde{\alpha}+\widetilde{\beta})^2}{\widetilde{\alpha}},
 \ \ \ \  (\widetilde{\alpha}:=\sqrt{\alpha^2+k\beta^2}, \ \
 \widetilde{\beta}:=\epsilon \beta).$$
It has been proven in \cite{SYa} that a square metric  is
Riemannian or flat-parallel if it is of constant flag curvature on
a closed manifold.

\bigskip

\noindent {\bf Case II:} In case of $n\ge 3$, suppose $F$ is of
scalar flag    curvature and $\beta$ is closed (or equivalently,
$F$ is locally projectively flat).

Under the deformation (\ref{gjcw61}), where
$u=u(b^2),v=v(b^2),w=w(b^2)$ satisfy the ODE
(\ref{gjcw58})--(\ref{gjcw60}), we have proved that $h$ is of
constant sectional curvature $\mu$ and $\rho$ satisfies
(\ref{gjcw63}) for some $c=c(x)$ and then (\ref{y9})
 holds.

  Let $h$ take the local form
 (\ref{gjcw67}).
 Now $c$ in (\ref{gjcw63})  can be expressed as
  \be\label{yy084}
 c=\frac{-k+\mu\langle a,x\rangle}{2\sqrt{1+\mu|x|^2}},
  \ee
  where $k,a$ are the same as that in (\ref{gjcw68}) and
  (\ref{gjcw69}).

{\bf Case IIA:} Assume $\mu<0$. We have $c=0$ by Lemma
\ref{lem61}. By (\ref{yy084}), we have $k=0,a=0$. Then by
(\ref{gjcw68}) we get $\rho=0$. Thus (\ref{gjcw61}) shows
$\beta=0$. Therefore, $F=\alpha$ is Riemannian in this case.

{\bf Case IIB:} Assume $\mu=0$. We have $c=-k/2=constant$. We will
show that $k=0$. Assume $k\ne 0$. Since $\alpha,\beta$ are defined
on the whole $M$, (\ref{gjcw61}) defines a Riemann metric $h$ and
a 1-form $\rho$ on $M$ for a suitable choice of $u,v,w$. For
example, we may choose $u,v,w$ as in (\ref{gjcw071}) here.
Therefore, the scalar function $||\rho||_h^2$ is defined on the
whole $M$. Then by (\ref{gjcw69}) and $\mu=0$, the scalar function
$f:=||\rho||_h^2$ is globally defined on $M$, and $f$ can be
locally expressed as
 $$f=k^2|x|^2+2k\langle a,x\rangle+|a|^2.$$
Since $\mu=0$, we have
 $$f_{|i|j}=f_{x^ix^j}=2k^2\delta_{ij}=2k^2h_{ij},$$
where the covariant derivatives are taken with respect to $h$. By
the above we have
$$\Delta f=2k^2n,$$
where $\Delta$ is the Laplacian of $h$. Integrating the above on
$M$ yields
 $$0=\int_M \Delta f \ dV_h=\int_M 2k^2n \ dV_h=2k^2n\ Vol_{h}(M).$$
Obviously it is a contradiction since $k\ne 0$. So we have $k=0$.
Then by (\ref{gjcw69}) we have $p^2=||\rho||_h^2$ is a constant on
$M$, and thus by (\ref{gjcw62}) we have $b^2=constant$ on $M$. Now
by (\ref{gjcw61}) we get
 \be\label{gjcw80}
 \alpha^2=\frac{1}{u}h^2-\frac{v}{uw^2}\rho^2,\ \
 \beta=\frac{1}{w}\rho.
 \ee
Thus (\ref{gjcw80}) shows that $\alpha$ is flat and $\beta$ is
parallel with respect to $\alpha$, since locally $\alpha$ and
$\beta$ are independent of $x\in M$, which follows from the fact
that $u=u(b^2),v=v(b^2),w=w(b^2)$ are constant and locally $h=|y|$
and $\rho=\langle a,y\rangle$.

{\bf Case IIC:} Assume $\mu>0$. By (\ref{gjcw68}) and
(\ref{yy084}), we can easily rewrite $\rho=p_0$ as
 \be\label{gjcw81}
\rho=\frac{2}{\mu}\ c_0,
 \ee
where $c_i:=c_{x^i},c_0:=c_iy^i$. Now by (\ref{gjcw80}) and
(\ref{gjcw81}) we easily get (\ref{ycw11}). Finally, it follows
from (\ref{gjcw63}) and then from (\ref{gjcw66})that $\tau$ and
$\lambda$ are given by (\ref{ycw12}).  \qed

\section{Proof of Theorem \ref{th3}}
In this section, we will prove Theorem \ref{th3}.  In the
following discussion,  we assume $\epsilon=1$. The case
$\epsilon=-1$ is similar.

 Note that
(\ref{y0010}) naturally holds if
$F=\alpha+a_1\beta/\alpha+\epsilon \beta^2/\alpha$ is a square
metric since in this case we have $a_1=\pm 2,\epsilon=1,k_1=2$.
Then the conclusion in Theorem \ref{th3}(ii) follows from Theorem
\ref{th2}(i) when $F$ is of constant flag curvature for $n\ge 3$.
In fact, when $F$ is of constant flag curvature for $n\ge 2$,
Theorem \ref{th3}(ii) follows directly  from \cite{SYa} if $F$ is
a square metric, and form \cite{CST} if $F$ is not a square
metric.

Since $\phi(s)=1+a_1s+\epsilon s^2$ is a solution of (\ref{y1}),
 by Theorem \ref{th1}, in the following proof we only need to
show that $\beta$ is closed when $F=\alpha+a_1\beta+\epsilon
\beta^2/\alpha$ is of scalar flag curvature in case of $n\ge 3$.

\begin{lem}\label{lem70}
 $\beta$ is closed $\Longleftrightarrow t_{ij}=0$ $\Longleftrightarrow t^k_{\
 k}=0$.
\end{lem}

\begin{lem}\cite{SYa}\label{lem71}
 For a scalar function $c=c(x)$, there does not hold
  $$\alpha b_k-sy_k\equiv 0 \ \ \ \  mod \ \  (s+c)$$
 for all $k$.
\end{lem}

Now we begin our discussion. We will prove our results in two
cases: $a_1\ne 0$ and $a_1=0$.

\bigskip
\

\noindent {\bf Case I:}  Assume $a_1\ne 0$.
\bigskip

We may assume $a_1>0$. Multiply $W^i_{\ k}=0$ by
 $$(n^2-1)(1-s^2)^4(1+2b^2-3s^2)^5\alpha^4$$
and we get an equation which is denoted by
 \beq\label{gjcw82}
 Eq_1:&=&648(n-2)(1+s)^4(1-s)^4s^3(\alpha
 b_k-sy_k)y^i\big[(s^2-1)r_{00}+2(a_1+2s)\alpha s_0\big]^2\nonumber\\
 &&+C^i_k(1+2b^2-3s^2)=0,
 \eeq
 where $C^i_k$ can be written in the form (\ref{ycw}).

\begin{lem}\label{lem72}
 Suppose
 \be\label{gjcw83}
  (s^2-1)r_{00}+2(a_1+2s)\alpha s_0 \equiv 0 \ \ \ \  mod \ \ (1+2b^2-3s^2).
  \ee
 Then we have
 \be\label{gjcw84}
  r_{00}=\tau\alpha^2(1+2b^2-3s^2), \ \ \ \ s_0=0,
  \ee
  where $\tau=\tau(x)$ is a scalar function.
\end{lem}

{\it Proof} : Eq. (\ref{gjcw84}) implies that there are
$f_0,f_1,f_2$ satisfying
 $$\alpha^2\big[(s^2-1)r_{00}+2(a_1+2s)\alpha s_0\big]=(f_0+f_1\alpha+f_2\alpha^2)\alpha^2(1+2b^2-3s^2),$$
which is equivalent to
 $$(2a_1s_0-f_1-2f_1b^2)\alpha^3+3f_1\beta^2\alpha=0,$$
 $$-(1+2b^2)f_2\alpha^4+(4\beta s_0-r_{00}-f_0-2b^2f_0+3f_2\beta^2)\alpha^2+\beta^2(3f_0+r_{00})=0.$$
Now by the above two equations, we can easily show that
(\ref{gjcw84}) holds. \qed

\bigskip

Now  we have (\ref{gjcw84}) by (\ref{gjcw82}), Lemma \ref{lem71}
and Lemma \ref{lem72}. Then by (\ref{gjcw84}) we can obtain the
expressions of the following quantities:
 $$
r_{00},\ r_i,\ r^m_m,\ r, \ r_{00|0}, \ r_{0|0}, \ r_{00|k}, \
r_{k0|m}, \ r_{k|0}, \ s_{k|m}, \ t_k,\ q_{km}, \ q_k,\ b^mq_{0m},
\ etc.
 $$
 For example, we have
 $$
 r_{00|0}=(1+2b^2-3s^2)\tau_0\alpha^2-2s(1+8b^2-9s^2)\tau^2\alpha^3,
 \ \ \tau_i:=\tau_{x^i},
 $$
and
 $$
 s_{k|m}=0, \ \ t_k=0, \ \ q_{km}=\tau (1+2b^2)s_{km},
 \ \ q_{00}=0, \ \ q_{k}=0, \ \ b^mq_{0m}=0.
$$
Plug all the above quantities into (\ref{gjcw82}) and then
multiplied by $1/(1+2b^2-3s^2)^5$ the equation (\ref{gjcw82}) can
be written as
 \be\label{gjcw85}
 Eq_2:=D^i_k(1-s)+12(n+1)(a_1+2)^2\alpha^2(\alpha b_k-y_k)y^it_{00}=0,
 \ee
 where $D^i_k$ can be written in the form (\ref{ycw}).
Similarly by (\ref{gjcw85}) we obtain

 \be\label{gjcw86}
  t_{00}=\gamma (\alpha^2-\beta^2),
  \ee
  where $\gamma=\gamma(x)$ is a scalar function.
Then by (\ref{gjcw86}) we have
 \be\label{gjcw87}
 t_{i0}=\gamma (y_i-\beta b_i), \ \ \ t^m_m=\gamma (n-b^2).
 \ee
Plug (\ref{gjcw86}) and (\ref{gjcw87}) into (\ref{gjcw85}) and
then multiplied by $1/(1-s)$ the equation (\ref{gjcw85}) is
written as
 \be\label{gjcw88}
 Eq_3=0.
 \ee
It follows from
 $$Eq_3 \equiv 0 \ \ \ \  mod \ (1-s)$$
that
 \be\label{y41}
 3(n-1)s_{i0}s_{k0}+\gamma (\alpha b_k-y_k)\big[(n-1)\alpha
 b_i-(3+b^2)y_i\big] \equiv 0 \ \ \ \  mod \ (1-s).
 \ee
The following lemma has been proved in \cite{SYa}. Now we have
proved that $\beta$ is closed.

\begin{lem}
 Suppose (\ref{y41}) holds for some scalar $\gamma=\gamma(x)$.  Then $\beta$ is closed.
\end{lem}

\bigskip

\

\noindent {\bf Case II:}  Assume $a_1= 0$. In this case, there are
some different steps from that in Case I.

\bigskip

Multiply $W^i_{\ k}=0$ by
 $$(n^2-1)(1-s^2)^4(1+2b^2-3s^2)^5\alpha^4$$
and  firstly we get (\ref{gjcw82}) with $a_1=0$, and then we have
(\ref{gjcw83}) with $a_1=0$. By a similar proof as that in Lemma
\ref{lem72} we obtain

\be\label{gjcw90}
  r_{00}=\tau\alpha^2(1+2b^2-3s^2)+\frac{6}{1-b^2}\beta s_0,
  \ee
  where $\tau=\tau(x)$ is a scalar function. Likely,  by (\ref{gjcw90}) we can obtain the
expressions of the following quantities:
 $$
r_{00},\ r_i,\ r^m_m,\ r, \ r_{00|0}, \ r_{0|0}, \ r_{00|k}, \
r_{k0|m}, \ r_{k|0}, \ s_{k|m}, \ t_k,\ q_{km}, \ q_k,\ b^mq_{0m},
\ etc.
 $$
Plug all the above quantities into (\ref{gjcw82}) with $a_1=0$,
and then multiplied by $(1-b^2)^3/(1+2b^2-3s^2)^5$ the equation
(\ref{gjcw82}) with $a_1=0$ can be written as
 \be\label{gjcw91}
 D^i_k(1-s)+48(n+1)(1-b^2)^2\alpha^2(\alpha b_k-y_k)y^i\big[(1-b^2)t_{00}-s_0^2\big]=0,
 \ee
where $D^i_k$ can be written in the form (\ref{ycw}). Then by
(\ref{gjcw91}) we have
 \be\label{gjcw92}
 (1-b^2)t_{00}-s_0^2 \equiv 0 \ \ \ \  mod \ (1-s).
 \ee
Now by (\ref{gjcw92}) we easily get
 \be\label{gjcw93}
 t_{00}=\gamma (\alpha^2-\beta^2)+\frac{s_0^2}{1-b^2},
 \ee
where $\gamma=\gamma(x)$ is a scalar function. It follows from
(\ref{gjcw93}) that
 \be\label{gjcw94}
 t_{k0}=\gamma (y_k-\beta b_k)+\frac{s_ks_0}{1-b^2}, \ \ \ t^m_m=\gamma
 (n-2b^2),\ \ t_k=(1-b^2)\gamma b_k, \ \ s_ms^m=-b^2(1-b^2)\gamma
 .
 \ee
Plug (\ref{gjcw93}) and (\ref{gjcw94}) into (\ref{gjcw91}) and
then multiplied by $1/(1-s^2)$ the equation (\ref{gjcw91}) is
written as
 \be\label{gjcw95}
 A^i_k(\alpha^2-\beta^2)+B^i_k=0,
 \ee
where $A^i_k$ and $B^i_k$ are polynomials in $(y^i)$. Contracting
(\ref{gjcw95}) by $b_ib^k$ we obtain
 \be\label{gjcw96}
 X\alpha^2-\beta^2\big[6(1-b^2)^3\gamma
 \beta^2-2(n^2-3n+5)(1-b^2)s_0^2+X\big]=0, \ \ (X:=-A^i_kb_ib^k),
 \ee
By (\ref{gjcw96}), there is a scalar function $\xi=\xi(x)$ such
that $X=\xi \beta^2$, and then plugging it into (\ref{gjcw96})
yields
 \beq
 \xi(\alpha^2-\beta^2)-6(1-b^2)^3\gamma
 \beta^2+2(n^2-3n+5)(1-b^2)s_0^2&=&0,\label{gjcw97}\\
\xi(a_{ij}-b_ib_j)-6(1-b^2)^3\gamma
 b_ib_j+2(n^2-3n+5)(1-b^2)s_is_j&=&0.\label{gjcw98}
 \eeq
Summing (\ref{gjcw98}) over the indices $i,j$ and using the
expression of $s_ms^m$ in (\ref{gjcw94}) we obtain
 \be\label{gjcw99}
 \xi(n-b^2)-2b^2((1-b^2)^2(n^2-3n+8-3b^2)\gamma=0.
 \ee
Contracting (\ref{gjcw98}) by $b^ib^j$ gives
 \be\label{gjcw100}
 \xi=6b^2(1-b^2)^2\gamma.
 \ee
Substitute (\ref{gjcw100}) into (\ref{gjcw99}) and we have
 \be\label{gjcw101}
 2(n-2)(n-4)b^2(1-b^2)^2\gamma=0.
 \ee

{\bf (1).} If $n\ne 4$, then by (\ref{gjcw101}) we have
$\gamma=0$. In this case, by $\gamma=0$ and (\ref{gjcw97}) we have
$s_0=0$ and then by (\ref{gjcw93}) we get $t_{00}=0$. Therefore,
it shows $\beta$ is closed by Lemma \ref{lem70}.

{\bf (2).} If $n= 4$, then plugging $n= 4$ and (\ref{gjcw100})
into (\ref{gjcw97}) shows
 \be\label{gjcw102}
 (1-b^2)\gamma (b^2\alpha^2-\beta^2)+3s_0^2=0.
 \ee
Clearly by (\ref{gjcw102}) we have $s_0=0,\gamma=0$, and then
again by (\ref{gjcw93}) we get $t_{00}=0$. Therefore, $\beta$ is
closed by Lemma \ref{lem70}.

\bigskip

Finally, we show (\ref{ycw15}) and (\ref{ycw16}). Since
$\phi(s)=1+a_1s+\epsilon s^2$, we have
$k_1=2\epsilon,k_2=0,k_3=-3\epsilon$ in
(\ref{gjcw58})--(\ref{gjcw60}). So in this case, we may choose
 \be\label{gjcw103}
u=(1-\epsilon b^2)^2,\ v=0, \ w=\sqrt{1-\epsilon b^2}.
 \ee
Then by (\ref{gjcw103}), (\ref{gjcw62}), (\ref{gjcw69}) and
(\ref{yy084}) we have
 \be\label{gjcw104}
p^2=\frac{4(\delta^2-\mu c^2)}{\mu^2},\ \
b^2=\frac{p^2}{1+\epsilon p^2}=\frac{4(\delta^2-\mu
c^2)}{\mu^2+4\epsilon(\delta^2-\mu c^2) },
 \ee
where $\delta$ is defined by (\ref{y010}) and  $\delta$  is a
constant by
 Lemma \ref{lem61}. Now by
(\ref{gjcw103}) and (\ref{gjcw104}) we obtain (\ref{ycw15})  from
(\ref{ycw11}), and using (\ref{gjcw103}), (\ref{gjcw104}) and
$k_1=2\epsilon,k_2=0,k_3=-3\epsilon$, we first get $\lambda,\tau$
by (\ref{ycw12}) and then we obtain (\ref{ycw16}) from (\ref{y6}).
\qed

\begin{rem}
Under the choice of $u,v,w$ as in (\ref{gjcw103}), we can evaluate
 $c$ and $\beta/\alpha$. Assume $\mu>0$ and we only show the
case $\epsilon=-1$. Note that
$F=\alpha+a_1\beta/\alpha-\beta^2/\alpha$ is a (regular) Finsler
metric if and only if $b^2<1/2$. So we have
 \be\label{gjcw105}
 \frac{2}{3}<1-p^2=\frac{1}{1+b^2}\le 1.
 \ee
 By (\ref{gjcw105}) we get
  \be\label{gjcw106}
\frac{2}{3}<1-p^2=\frac{\mu^2-4\delta^2+4\mu c^2}{\mu^2}\le 1
  \ee
  Then (\ref{gjcw106}) shows
   \be\label{gjcw107}
 \frac{\delta^2}{\mu}-\frac{\mu}{12}<c^2\le \frac{\delta^2}{\mu}.
   \ee
   Now by (\ref{ycw15}) we obtain
    \be
  |\frac{\beta}{\alpha}|\le \frac{2\sqrt{\delta^2-\mu c^2}}{\mu^2-4\delta^2+4\mu
  c^2}< \frac{1}{\sqrt{2}}.
    \ee
\end{rem}

\begin{rem}\label{rem72}
For the parameters $k_1,k_2,k_3$ in (\ref{y1}) satisfying certain
equations, $F=\alpha\phi(\beta/\alpha)$ can be essentially
considered as the metric in Theorem \ref{th3}. For example, if
 $k_1,k_2,k_3$ satisfy (\ref{gjcw76}), then we get (\ref{F}) which
 can be considered as a square metric if $k>-1/b^2$ in (\ref{F}).
 If
   $k_2=(2k_1+3k_3)(3k_1+2k_3)/25,$
   then $F=\alpha\phi(\beta/\alpha)$ can be a regular
   $(\alpha,\beta)$-metric in the form $F=\widetilde{\alpha}\pm
   \widetilde{\beta}^2/\widetilde{\alpha}$, where we define
   $\widetilde{\alpha}:=\sqrt{\alpha^2+k\beta^2},\widetilde{\beta}:=\epsilon\beta$
   with $k,\epsilon$ being constant satisfying $k>-1/b^2$. A
   two-dimensional $(\alpha,\beta)$-metric
   $F=\alpha\pm \beta^2/\alpha$ can be locally projectively flat with
   $\beta$ being not closed (\cite{Y1}).
\end{rem}

\bigskip

{\bf Acknowledgement:}

The  author expresses his sincere thanks to China Scholarship
Council for its funding support. This job was done during the
period (June 2012--June 2013) when he as a postdoctoral
researcher visited Indiana University-Purdue University
Indianapolis, USA.

\vspace{0.5cm}

\noindent Guojun Yang \\
Department of Mathematics \\
Sichuan University \\
Chengdu 610064, P. R. China \\
 ygjsl2000@yahoo.com.cn

\end{document}